\documentclass[11pt,letterpaper,nobysame]{amsart}

\usepackage{amsbsy,amssymb,amscd,amsfonts,latexsym,amstext,delarray,amsmath,mathtools}
\usepackage{mathabx} 
\usepackage{euscript} 
\usepackage{float}
\usepackage[alphabetic]{amsrefs}

\usepackage{graphicx}
\usepackage{tikz}
\usetikzlibrary{cd, matrix, arrows, positioning, arrows.meta, shapes.callouts, trees, backgrounds, fit, shapes.geometric, shapes.arrows, calc}
\usepackage{etex, pictexwd, dcpic} 
\input xy
\xyoption{all}

\usepackage[T1]{fontenc}
\usepackage{lmodern} 
\usepackage{microtype} 
\usepackage{charter} 
\usepackage{setspace}
\usepackage{fancyhdr}
\usepackage{parskip} 
\usepackage{enumitem} 
\usepackage{booktabs} 
\usepackage[normalem]{ulem} 
\usepackage{comment}
\usepackage{caption}
\usepackage{blkarray}
\usepackage{multirow}
\usepackage{array}
\usepackage{footnote}

\usepackage{hyperref}

\input xy
\xyoption{all}
\usepackage{hyperref}
\hypersetup{
    colorlinks=true,
    linkcolor=blue,
    filecolor=magenta,      
    urlcolor=cyan,
}

\theoremstyle{definition}

\theoremstyle{remark}

\def\R{\mathbb{R}}     
\def\C{\mathbb{C}}     

\begin{document}

{{

\title{\bf{{ADVANCING  MATHEMATICS RESEARCH\\ WITH GENERATIVE AI}}}

\date{\today}

\begin{abstract} 
The main drawback of using generative AI models for advanced mathematics is that these models are not primarily logical reasoning engines. However,  Large Language Models, and their refinements,  can pick up on patterns in higher mathematics that are difficult for humans to see. By putting the design of generative AI models to their advantage, mathematicians may use them as powerful interactive assistants that can carry out laborious tasks, generate and debug code, check examples, formulate conjectures and more. We discuss how generative AI models can be used to advance mathematics research. We also discuss their integration  with neuro-symbolic solvers, Computer Algebra Systems and  formal proof assistants such as Lean.
\end{abstract}

\author[Lisa Carbone]{Lisa Carbone}
\address{Department of Mathematics, Rutgers University, Piscataway, NJ 08854-8019, USA}
\email{lisa.carbone@rutgers.edu}

\maketitle


\section{Introduction}

Mathematicians have mixed views about the role of generative AI models in the mathematical landscape. While these models, such as Large Language Models (LLMs), can look convincingly like they replicate known mathematics, on careful scrutiny, it is apparent that they are just masters of the {\it rhetoric} of mathematics. They don't meet  the standards of rigor that the field requires.

This limitation is inherent to the architecture of LLMs. They  are fundamentally statistical, not logical, engines. As next-word predictors, they lack a built-in engine for formal symbolic deduction, meaning their emergent `reasoning', called {\it Natural Language Reasoning}, is derived  from statistical patterns in their training data.

However, this statistical foundation also provides certain benefits,  as LLMs can pick up on patterns in higher mathematics that are difficult for humans to see. In particular, they have a learned geometric representation of mathematical language. By leveraging these capabilities and working with LLMs as they were designed to operate, mathematicians can use them as powerful interactive assistants. 

The AI landscape is rapidly evolving beyond standard LLMs, with Large Reasoning Models (LRMs) and Large Context Models (LCMs) emerging as the next generation of this technology. 
While built on the foundation of LLMs, these models have distinct architectures and methods of training.
In particular, LRMs, such as Gemini,  use hybrid neuro-symbolic systems which combine text generation with verification (see Subsection~\ref{LRM-LCM}).

 In Section~\ref{GenAI}, we describe the  capabilities and limitations of LLMs and  the nature of their training data.
  In Section~\ref{mechanics}, we explore the mechanics of generative AI, from the management of context windows to the geometric representation of mathematical language within high-dimensional embedding spaces.   We also discuss the evolution of architecture from standard language models to Large Reasoning Models (LRMs) that utilize neuro-symbolic verification.

In Sections~\ref{capabilities} through ~\ref{output}, we examine methods for influencing the outputs of  these models. We discuss how they simulate computational environments and how prompt engineering can impose logical order on probabilistic outputs (Section~\ref{prompts}).

In Section~\ref{CGT}, we discuss  how generative AI models are useful for questions in  Combinatorial Group Theory.
In Section~\ref{hybrid}  we discuss the integration of LLMs with Computer Algebra Systems (CAS) and formal proof assistants such as Lean.

Finally, in Section 9, we imagine a collaborative research partnership where a  mathematician engages in an iterative dialogue with AI to discover new ideas and research directions.

\section{Generative AI models}\label{GenAI}
\subsection{Mathematical training data for generative AI models}

Training data for generative AI models includes web documents,  code with mathematical content, textbooks, online faculty-authored lecture notes and course materials, solutions to problem sets, as well as   academic and scientific journal papers and content from arXiv. Google's data collection from Google Books, for example, is unparalleled in its scale.

The training data includes the  standard undergraduate mathematics curriculum from US and other universities, as well as the standard coursework curriculum for PhD programs in mathematics from US and other universities.

Differences in the content and handling of training data of AI companies are primary reasons why various AI models have distinct strengths, weaknesses, and `personalities', especially in a specialized domain like mathematics. Curation and filtering of data, such as removing low-quality content and emphasizing  trusted sources, is essential.  However, the training material for each model and its handling remains a closely guarded trade secret.

A key limitation  is that a generative AI model's knowledge base is heavily biased towards materials that are easily digitized.

\subsection{How much do generative AI models `know'?}
Current generative AI models are trained on human data, which means they cannot generate knowledge outside of that data. While they can go beyond the knowledge of any single human, they can only come up with new ideas through extrapolation, not discovery from first principles. Their core capability is {\it Natural Language Reasoning} - the ability of generative AI models to process, synthesize and generate human language by identifying and replicating statistical patterns from vast amounts of data in the form of tokens.  

In contrast,  {\it Symbolic Deduction} is  manipulation based on formal rules applied to abstract symbols, as in the operation of a Computer Algebra System. Current LLMs are built to imitate reasoning based on patterns, but are incapable of performing true symbolic deduction.

However, just as mathematicians are becoming familiar with the use of generative AI models,  their internal architecture is already changing.

In 2024, Google's DeepMind reported silver-medal level on  International Mathematical Olympiad problems using AlphaProof and AlphaGeometry. In 2025, DeepMind reported gold-medal level using Gemini `Deep Think'.  These systems  did not use a separate formal proof assistant.  The key innovation in their new model is the integration of Symbolic Deduction with Natural Language Reasoning, representing a first step towards an AI-based alternative to formal proof assistants.  Their new models  also incorporate  techniques like parallel thinking, which allowed them to explore multiple solutions simultaneously.

The significant advancement in reasoning capability has been achieved by neuro-symbolic systems that use an advanced form of Chain-of-Thought and Tree-of-Thought reasoning with a feedback loop for verification.  

The next frontier goes even further, with systems like Google's AlphaZero, which is not based on human-generated data. AlphaZero is given only the rules  and  foundations of a subject. It then learns (for example, in chess) by self-play and  it generates its own training data. It is therefore likely that a future `AlphaMath' system could be trained only on the axioms and rules of inference of ZFC set theory. Such a system would engage in genuine `self-discovery' by attempting to prove theorems. Over time, with human feedback, it would begin to recognize patterns and strategies that lead to successful proofs.

\subsection{Drawbacks of using LLMs in mathematics}
\label{drawbacks}

The primary drawback of using LLMs for mathematics is that they are probabilistic pattern-matchers, not logical reasoning engines. Their core function is to predict the next most likely word or symbol in a sequence, not to apply deterministic  mathematical rules.  Generative AI models are  experts  at reconstructing what a proof should `look like'. This  limitation leads to several key problems:

{\bf Mathematical hallucinations:} AI models generate outputs that look plausible and are often formatted correctly but contain nonsensical logic, invented theorems, or critical errors. They can state these falsehoods with the same confident authority as factual information.

{\bf Lack of true reasoning:} Unlike a Computer Algebra System, an LLM cannot manipulate symbolic expressions according to deterministic rules. This makes it prone to subtle but serious errors in algebra, analysis, arithmetic, and logic. It is fundamentally incapable of handling verification.

{\bf Propagation of errors:}  AI models learn from their  training data, which includes errors and misconceptions found online. For example, it is known that ChatGPT 4.0 contains training data from retracted scientific papers \cite{MITTR2025}. An LLM will inherit and reproduce these mistakes, making them unreliable for tasks that require rigor.

\section{Foundations of Generative AI}\label{mechanics}

\subsection{Tokens: the building blocks of language}
In AI, a {\it token} is the fundamental building block of language. AI models convert input materials (prompts) into tokens. A token $t_i$  is a  digital representation of a word, a symbol,  punctuation, a diagram, or other input.

Let $V$ be a finite vocabulary consisting of token representations of the training data for an LLM. An {\it input} to an LLM is a sequence of tokens $(t_1, t_2, \dots, t_n)$. An {\it output} is a sequence of tokens $(t_1, t_2, \dots, t_m)$, where each $t_i \in V$.

The fundamental operation of an LLM is to compute the probability of the next token, $t_{k+1}$, conditioned on the sequence of all preceding tokens $(t_1, t_2, \dots, t_k)$.

All AI language models have a maximum token limit for input and output. The token input limit has an impact on the work possible in a single chat session.

\subsection{The LLM as a probabilistic knowledge graph}

Consider the following simplified  view of a generative AI model as a `knowledge graph' which represents all possible connections and probabilities within the LLM's knowledge base.

Each word is a node in the graph. Each edge is a statistical relationship or an association between words, which cluster to form concepts.

\bigskip

\begin{center}
\begin{figure}[H]
\begin{tikzpicture}[
    concept/.style={rectangle, rounded corners, draw, thick, fill=blue!10, minimum height=1cm, minimum width=2.2cm, font=\sffamily},
    active_prompt/.style={concept, fill=yellow!40, very thick, draw=orange},
    output_node/.style={circle, draw, thick, fill=purple!10, minimum size=1.5cm, font=\sffamily},
    high_prob/.style={-Latex, thick, draw=red!80},
    mid_prob/.style={-Latex, thick, draw=gray!80},
    low_prob/.style={-Latex, dashed, draw=black!80, thin},
    desc/.style={align=left, font=\sffamily\small, text width=9cm}
]
\node[active_prompt] (concept_a) at (0,2) {addition};
\node[active_prompt] (concept_b) at (6,2) {negative numbers};
\node[concept] (concept_c) at (12,2) {irrational numbers};
\node[output_node, fill=green!20, very thick, draw=green!60] (node_x) at (0.5,-1.3) {integers};
\node[output_node] (node_y) at (7.5,-1.3) {real numbers};
\node[output_node] (node_z) at (3,-3.9) {inverse};
\draw[-{Stealth[length=5mm]}, line width=2.5pt, draw=red!90] (concept_a) to[bend left=20] (node_x);
\draw[-{Stealth[length=5mm]}, line width=2.5pt, draw=red!90] (concept_b) to (node_x);
\node[below=2pt of node_x, red!90, font=\sffamily\bfseries] {{\tiny{CONVERGED PROBABILITY}}};
\draw[low_prob] (concept_a) to[bend right=20] (node_y);
\draw[low_prob] (concept_c) to (node_y);
\draw[low_prob] (concept_b) to (node_z);
\end{tikzpicture}
\caption{A simplified  view of an LLM as a `knowledge graph'. }
\end{figure}
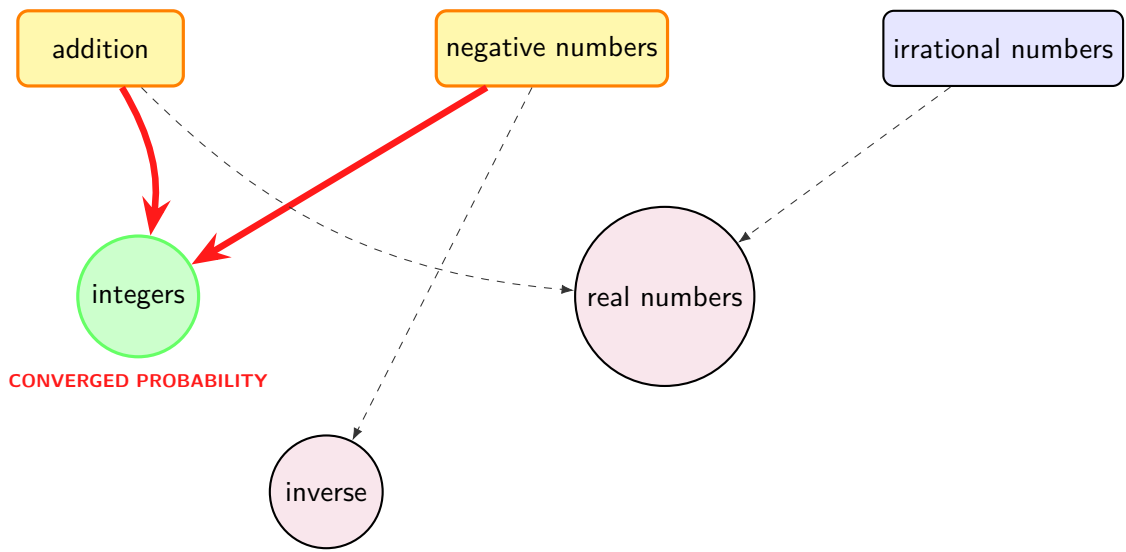
\end{center}

\bigskip
Weighted edges between trillions of nodes represent how words or concepts are related. The weights are probabilities indicating the likelihood of  connection. 
 
An LLM generates text by predicting the next word, starting from the prompt, while calculating the probabilities of all possible next nodes. It then chooses one of the most likely paths through the knowledge graph.

\subsection{Generative AI models and patterns in higher mathematics}
\label{patterns}

A generative AI model can identify mathematical patterns that are not obvious to humans.  It performs massive scale statistical analysis on the symbolic structure of mathematics itself. This allows it to simulate a mathematician's intuition.  It can find cross-disciplinary correlations and potential analogies that humans may not think of.

The way mathematicians are trained to think assumes the progression of ideas of mathematical discovery. We have learned concepts in a chronological sequence.   AI models do not `know' that mathematical fields developed sequentially. For an AI model, mathematical concepts coexist simultaneously. They are linked only by the statistical and structural patterns embedded in the language of mathematics. 
Due to their inherent design, they can propose unexpected pathways, connections and points of view.

There are advantages of this capability. A generative AI model sees mathematics as a language with a complex grammar and vocabulary composed of  symbols and diagrams. It learns the statistical relationships between the tokens in this language.

Trained on the entire available library of mathematics and science without disciplinary boundaries, an LLM learns the statistical relationships between tokens, allowing it to explore the entire space of mathematical connections, not just those that individual mathematicians are trained to look for.

\newpage
\subsection{The embedding space for an LLM}

\begin{center}
\begin{figure}[H]
\begin{tikzpicture}
[
    font=\small,
    dot/.style={fill, circle, inner sep=1.5pt},
    map_arrow/.style={->, gray, thin, -{Stealth[length=6pt]}},
    vec_arrow/.style={->, red, thick, dashed, -{Stealth[length=8pt]}}
]
    \node (V_box) [
        draw,
        circle,
        minimum height=4cm,
        minimum width=4cm,
        label={[font=\large]above:Vocabulary $V$}
    ] at (0,0) {};

    \node (man) at (V_box.center) [yshift=-0.5cm] {man};
    \node (woman) [below=0.1cm of man] {woman};
    \node (queen) [above=0.1cm of man] {queen};
    \node (king) [above=0.1cm of queen] {king};

    \begin{scope}[xshift=8cm, yshift=-2.5cm] 
        \node [font=\large, align=center] at (2.5, 5.5) {Embedding Space $\mathbb{R}^d$};

        \draw [->, thick] (-0.5, 0) -- (5, 0) node[right] {}; 
        \draw [->, thick] (0, -0.5) -- (0, 5) node[above] {}; 


        \node (man_vec)   [dot] at (1, 1) {};
        \node (king_vec)  [dot] at (2, 3) {};
        \node (woman_vec) [dot] at (3, 1.5) {};
        \node (queen_vec) [dot] at (4, 3.5) {};

        \node [below=0.0pt of man_vec] {$\varphi(\text{man})$};
        \node [above left=2pt of king_vec] {$\varphi(\text{king})$};
        \node [below right=0.3pt of woman_vec] {$\varphi(\text{woman})$};
        \node [above=6pt of queen_vec] {$\varphi(\text{queen})$};

        \draw [vec_arrow] (man_vec) -- (king_vec)
            node [midway, right, text=black, xshift=2mm] {$\vec{v}_{\text{royal}}$};
        \draw [vec_arrow] (woman_vec) -- (queen_vec)
            node [midway, right, text=black, xshift=2mm] {$\vec{v}_{\text{royal}}$};

        \draw [->, blue, thick, dashed, -{Stealth[length=8pt]}] (man_vec) -- (woman_vec)
            node [midway, below, text=black, yshift=-0.4mm] {$\vec{v}_{\text{gender}}$};
        \draw [->, blue, thick, dashed, -{Stealth[length=8pt]}] (king_vec) -- (queen_vec)
            node [midway, above, text=black, yshift=0.4mm] {$\vec{v}_{\text{gender}}$};
    \end{scope}

    \node[font=\large] at (4.5, 2.5) {$\varphi: V \to \mathbb{R}^d$}; 

    \draw [map_arrow] (king.east) .. controls (2, 1.5) and (6.5, 3.5) .. (king_vec.west);
    \draw [map_arrow] (queen.east) .. controls (2, 0) and (6.5, 4) .. (queen_vec.west);
    \draw [map_arrow] (man.east) .. controls (2, -1) and (6.5, 0.0) .. (man_vec.west);
    \draw [map_arrow] (woman.east) .. controls (2, -2) and (6.5, 0.5) .. (woman_vec.west);
\end{tikzpicture}
\caption{A hypothetical embedding for $d=2$.}
\end{figure}
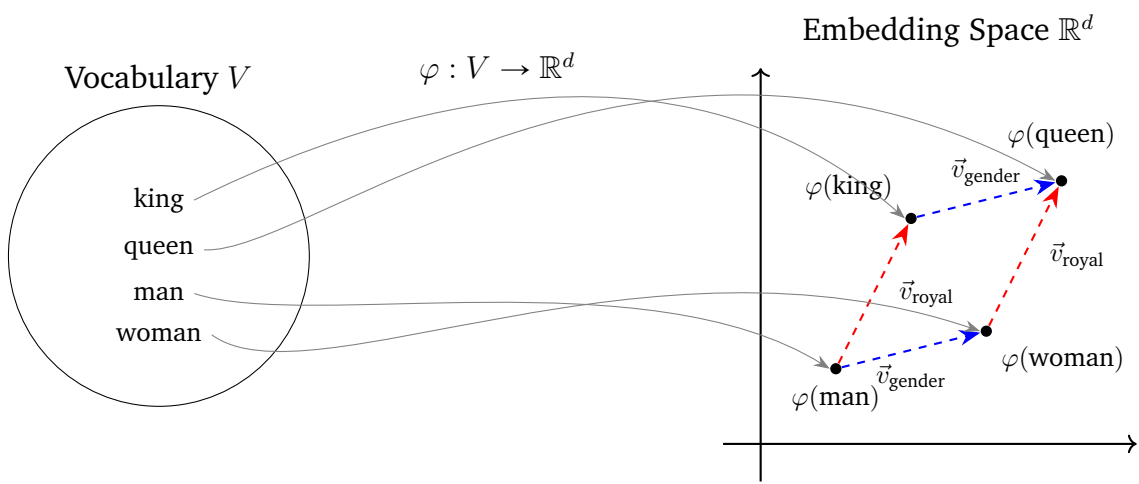
\end{center}

Generative AI models don't `see' symbols and equations. There is a map $\varphi:V\to \R^{d}$ (for $d\leq 13,000$ for large models) from the vocabulary $V$ of tokens to a high-dimensional vector space (the {\it embedding space}). 
The dimension $d$ is a model parameter. Relationships in the vocabulary  are translated into geometric relationships in~$\R^{d}$ \cite{mikolov2013}.

The model can detect geometric relationships  in $\R^{d}$ that are difficult for humans to visualize for $d>>0$ and that are not  obvious from the symbolic definitions. 

\subsection{The path of a word $w$ through an LLM}
Let $\varphi(w)$ denote the initial embedding of the token word $w$ in $\R^d$.   Its `meaning' is its relationship to all other images of tokens in the space and their locations. Furthermore, in advanced contextual models, this meaning changes with context.

The model uses  linear transformations  $T_q$, $T_k$ and $T_v$ to project each token embedding $\varphi(w)\in\mathbb R^d$ into three  specialized vectors which  have  different purposes: a query vector $q$ representing `what I am looking for', a key vector $k$ encoding `what topic I represent' and a value vector $v$ representing `what information I contain'. The maps
\[
T_q,\quad T_k,\quad T_v : \mathbb{R}^d \to \mathbb{R}^d,
\]
project $\varphi(w)$ to three specialized vectors:
\[
q = T_q \varphi(w), \qquad k = T_k \varphi(w), \qquad v = T_v \varphi(w).
\]

The model determines relevancy scores by measuring the dot product between the query vector $q$ of a token and the key vectors $k_i$ of every token in the sequence.
If $q$ and some $k_i$ point in similar directions, they obtain a high relevancy score.

It then uses these scores to create a new, context-rich vector as a weighted average of the value vectors. This new vector is then fed through a piecewise linear map to produce the  output vector, $\zeta^{(\ell)}(w)$ for level $\ell$ of the model.

 The process is repeated through many layers of the model;  the output $\zeta^{(\ell)}(w)$ becomes the input for layer $\ell+1$. The final contextualized representation, $\zeta^{(N)}(w)$ is the output of the final layer $N$ and is different for every sentence it appears in.

 The model generates an output sentence token-by-token in a
  loop. To predict the {next} token in the output, it uses the entire output sequence generated so far as its input. It then takes the final contextualized representation $\zeta^{(N)}(w_{s})$ of the {last} token $w_s$ from that sequence, passes it through one final linear map, and {adds} the most probable next token to the output sequence. This becomes the input for the  next step.

\begin{figure}[H]
\begin{tikzpicture}[
    scale=0.80,
    token_embed/.style={rectangle, draw, thick, minimum width=6cm, minimum height=1cm, align=center, fill=blue!10, font=\large},
    vector_box/.style={rectangle, draw, thick, minimum width=3cm, minimum height=1cm, align=center, fill=red!10, font=\large},
    label_style/.style={font=\small\itshape, text width=3.3cm, align=center, below=0.1cm},
    arrow_style/.style={-Latex, thick, rounded corners=5pt},
    layer_output/.style={rectangle, draw, thick, minimum width=6cm, minimum height=1cm, align=center, fill=orange!10, font=\large},
    transform_label/.style={font=\small, fill=white, inner sep=1pt, align=center},
    nn_node/.style={circle, draw, fill=white, inner sep=1.4pt},
    embed_label/.style={font=\small\itshape, fill=white, inner sep=1pt, align=center}
]

    \node (token) [token_embed] {Token Embedding $\varphi(w) \in \mathbb{R}^d$};

    \node[draw, dashed, rounded corners, thick, fit=(token), inner sep=6pt] (embedbox) {};
    \node[embed_label, right=0.25cm of embedbox.east] 
        {Embedding space $\mathbb{R}^d$};

    \node (query_vec) [vector_box, below left=2cm and 0.5cm of token.south west] {$q$};
    \node [label_style] at (query_vec.south) {Query: `what I am looking for'};

    \node (key_vec) [vector_box, below=2cm of token.south] {$k$};
    \node [label_style] at (key_vec.south) {Key: `what topic I represent'};

    \node (value_vec) [vector_box, below right=2cm and 0.5cm of token.south east] {$v$};
    \node [label_style] at (value_vec.south) {Value: `what information I contain'};

    \draw [arrow_style] (token.south) -- 
        node[transform_label, pos=0.5, left=2mm] {$T_q$\\\scriptsize linear transformation} 
        (query_vec.north);

    \draw [arrow_style] (token.south) -- 
        node[transform_label, pos=0.5, above=2mm] {$T_k$\\\scriptsize linear transformation} 
        (key_vec.north);

    \draw [arrow_style] (token.south) -- 
        node[transform_label, pos=0.5, right=2mm] {$T_v$\\\scriptsize linear transformation} 
        (value_vec.north);

    \node (layer_output) [layer_output, below=3.8cm of key_vec] 
        {Layer $\ell$ output $\zeta^{(\ell)}(w)$};

    \coordinate (stackcenter) at ($(layer_output.south)+(0,-3.8cm)$);

    \begin{scope}[shift={(stackcenter)}]

        \foreach \y [count=\i] in {1.8, 0.0, -1.8} {

            \node[nn_node] (in\i1) at (-2.4,\y+0.8) {};
            \node[nn_node] (in\i2) at (-2.4,\y) {};
            \node[nn_node] (in\i3) at (-2.4,\y-0.8) {};

            \node[nn_node] (hA\i1) at (-1.0,\y+1.1) {};
            \node[nn_node] (hA\i2) at (-1.0,\y+0.3) {};
            \node[nn_node] (hA\i3) at (-1.0,\y-0.3) {};
            \node[nn_node] (hA\i4) at (-1.0,\y-1.1) {};

            \node[nn_node] (hB\i1) at (0.9,\y+0.7) {};
            \node[nn_node] (hB\i2) at (0.9,\y) {};
            \node[nn_node] (hB\i3) at (0.9,\y-0.7) {};

            \node[nn_node] (out\i)  at (2.4,\y) {};

            \foreach \a in {1,2,3}{
                \foreach \b in {1,2,3,4}{
                    \draw (in\i\a) -- (hA\i\b);
                }
            }

            \foreach \a in {1,2,3,4}{
                \foreach \b in {1,2,3}{
                    \draw (hA\i\a) -- (hB\i\b);
                }
            }

            \foreach \b in {1,2,3}{
                \draw (hB\i\b) -- (out\i);
            }

            \draw (in\i1) .. controls (-1.7,\y+2.1) and (0.0,\y+1.9) .. (hB\i3);
            \draw (in\i3) .. controls (-1.7,\y-2.1) and (0.0,\y-1.9) .. (hB\i1);
        }

        \node[draw, rounded corners, thick,
              minimum width=7cm,
              minimum height=6.5cm] (stackbox) at (0,0) {};

    \end{scope}

    \node[label_style, right=0.7cm of stackbox, text width=5.5cm, align=left] 
        {Stack of transformer blocks: \\ (attention + feed-forward \\neural networks) \\ piecewise-linear maps,\\ repeated $N$ times};

    \node (final_output) [token_embed, below=2.0cm of stackbox] 
        {Final contextualized representation $\zeta^{(N)}(w)$};

    \draw[arrow_style] (layer_output.south) -- (stackbox.north);
    \draw[arrow_style] (stackbox.south) -- (final_output.north);

    \coordinate (layer_input_point) at (layer_output.north);
    \draw [arrow_style] (query_vec.south) |- (layer_input_point);
    \draw [arrow_style] (key_vec.south) -- (layer_input_point);
    \draw [arrow_style] (value_vec.south) |- (layer_input_point);

\end{tikzpicture}

\caption{Projection of a token embedding $\varphi(w)$ in the embedding space 
$\mathbb{R}^d$ into specialized vectors ($q,k,v$) via linear transformations 
$T_q, T_k, T_v$, generation of the layer-$\ell$ output $\zeta^{(\ell)}(w)$, 
and repeated processing through a stack of transformer (neural network) blocks, 
including the feed-forward piecewise-linear map, producing the final contextualized representation 
$\zeta^{(N)}(w)$.}
\end{figure}
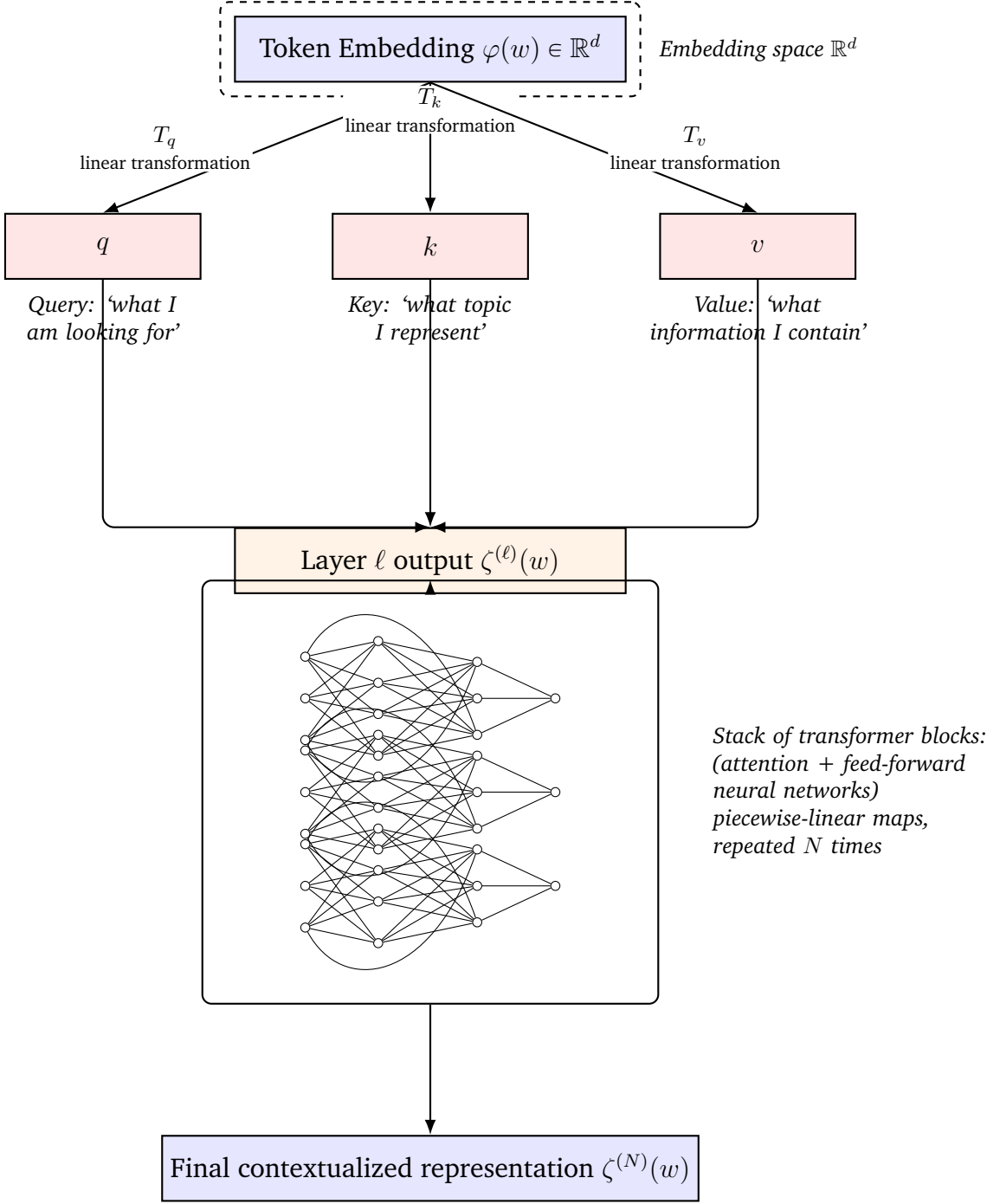

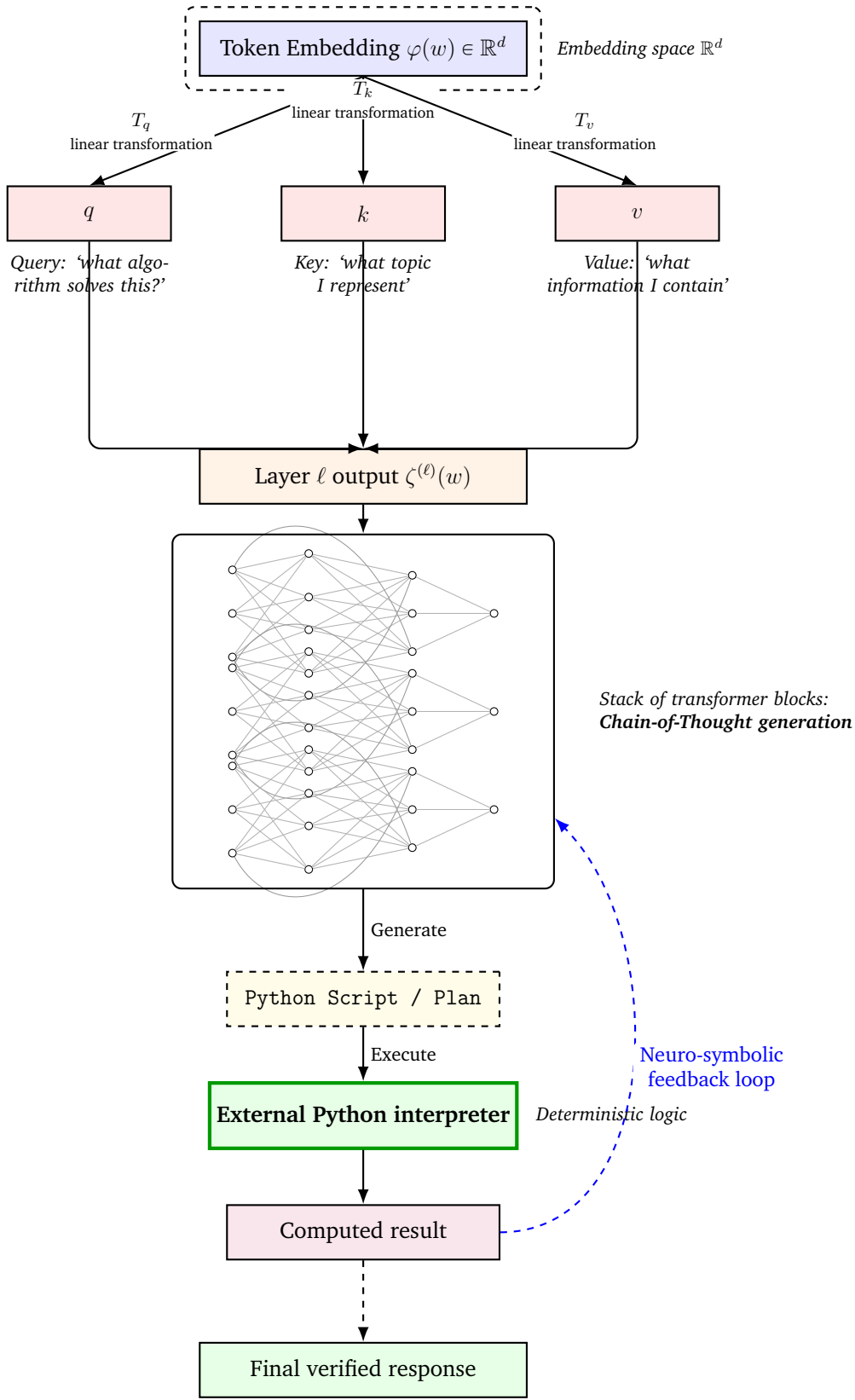
\begin{figure}[H]
\begin{tikzpicture}[
    scale=0.85, transform shape,
    token_embed/.style={rectangle, draw, thick, minimum width=6cm, minimum height=1cm, align=center, fill=blue!10, font=\large},
    vector_box/.style={rectangle, draw, thick, minimum width=3cm, minimum height=1cm, align=center, fill=red!10, font=\large},
    label_style/.style={font=\small\itshape, text width=3.3cm, align=center, below=0.1cm},
    arrow_style/.style={-Latex, thick, rounded corners=5pt},
    layer_output/.style={rectangle, draw, thick, minimum width=6cm, minimum height=1cm, align=center, fill=orange!10, font=\large},
    transform_label/.style={font=\small, fill=white, inner sep=1pt, align=center},
    nn_node/.style={circle, draw, fill=white, inner sep=1.4pt},
    embed_label/.style={font=\small\itshape, fill=white, inner sep=1pt, align=center},
    code_box/.style={rectangle, draw, thick, dashed, minimum width=5cm, minimum height=1cm, align=center, fill=yellow!10, font=\large\ttfamily},
    interpreter_box/.style={rectangle, draw=green!60!black, ultra thick, minimum width=5.5cm, minimum height=1.2cm, align=center, fill=green!10, font=\large\bfseries},
    result_box/.style={rectangle, draw, thick, minimum width=5cm, minimum height=1cm, align=center, fill=purple!10, font=\large}
]

    \node (token) [token_embed] {Token Embedding $\varphi(w) \in \mathbb{R}^d$};

    \node[draw, dashed, rounded corners, thick, fit=(token), inner sep=6pt] (embedbox) {};
    \node[embed_label, right=0.25cm of embedbox.east] 
        {Embedding space $\mathbb{R}^d$};

    \node (query_vec) [vector_box, below left=2cm and 0.5cm of token.south west] {$q$};
    \node [label_style] at (query_vec.south) {Query: `what algorithm solves this?'}; 

    \node (key_vec) [vector_box, below=2cm of token.south] {$k$};
    \node [label_style] at (key_vec.south) {Key: `what topic I represent'};

    \node (value_vec) [vector_box, below right=2cm and 0.5cm of token.south east] {$v$};
    \node [label_style] at (value_vec.south) {Value: `what information I contain'};

    \draw [arrow_style] (token.south) -- 
        node[transform_label, pos=0.5, left=2mm] {$T_q$\\\scriptsize linear transformation} 
        (query_vec.north);

    \draw [arrow_style] (token.south) -- 
        node[transform_label, pos=0.5, above=2mm] {$T_k$\\\scriptsize linear transformation} 
        (key_vec.north);

    \draw [arrow_style] (token.south) -- 
        node[transform_label, pos=0.5, right=2mm] {$T_v$\\\scriptsize linear transformation} 
        (value_vec.north);

    \node (layer_output) [layer_output, below=3.8cm of key_vec] 
        {Layer $\ell$ output $\zeta^{(\ell)}(w)$};

    \coordinate (stackcenter) at ($(layer_output.south)+(0,-3.8cm)$);

    \begin{scope}[shift={(stackcenter)}]
        \foreach \y [count=\i] in {1.8, 0.0, -1.8} {
            \node[nn_node] (in\i1) at (-2.4,\y+0.8) {};
            \node[nn_node] (in\i2) at (-2.4,\y) {};
            \node[nn_node] (in\i3) at (-2.4,\y-0.8) {};
            \node[nn_node] (hA\i1) at (-1.0,\y+1.1) {};
            \node[nn_node] (hA\i2) at (-1.0,\y+0.3) {};
            \node[nn_node] (hA\i3) at (-1.0,\y-0.3) {};
            \node[nn_node] (hA\i4) at (-1.0,\y-1.1) {};
            \node[nn_node] (hB\i1) at (0.9,\y+0.7) {};
            \node[nn_node] (hB\i2) at (0.9,\y) {};
            \node[nn_node] (hB\i3) at (0.9,\y-0.7) {};
            \node[nn_node] (out\i)  at (2.4,\y) {};
            \foreach \a in {1,2,3}{ \foreach \b in {1,2,3,4}{ \draw[black!30] (in\i\a) -- (hA\i\b); } }
            \foreach \a in {1,2,3,4}{ \foreach \b in {1,2,3}{ \draw[black!30] (hA\i\a) -- (hB\i\b); } }
            \foreach \b in {1,2,3}{ \draw[black!30] (hB\i\b) -- (out\i); }
            \draw[black!40] (in\i1) .. controls (-1.7,\y+2.1) and (0.0,\y+1.9) .. (hB\i3);
            \draw[black!40] (in\i3) .. controls (-1.7,\y-2.1) and (0.0,\y-1.9) .. (hB\i1);
        }
        \node[draw, rounded corners, thick,
              minimum width=7cm,
              minimum height=6.5cm] (stackbox) at (0,0) {};
    \end{scope}

    \node[label_style, right=0.7cm of stackbox, text width=5.5cm, align=left] 
        {Stack of transformer blocks: \\ \textbf{Chain-of-Thought generation} };

    \coordinate (layer_input_point) at (layer_output.north);
    \draw [arrow_style] (query_vec.south) |- (layer_input_point);
    \draw [arrow_style] (key_vec.south) -- (layer_input_point);
    \draw [arrow_style] (value_vec.south) |- (layer_input_point);
    \draw[arrow_style] (layer_output.south) -- (stackbox.north);


    \node (code) [code_box, below=1.5cm of stackbox] {Python Script / Plan};
    \draw [arrow_style] (stackbox.south) -- node[right, font=\small] {Generate} (code.north);

    \node (interpreter) [interpreter_box, below=1.0cm of code] {External Python interpreter};
    \draw [arrow_style] (code.south) -- node[right, font=\small] {Execute} (interpreter.north);
    \node [right=0.2cm of interpreter, font=\small\itshape, text width=3cm] {Deterministic logic };

    \node (result) [result_box, below=1.0cm of interpreter] {Computed result};
    \draw [arrow_style] (interpreter.south) -- (result.north);

    \node (final_output) [token_embed, below=1.5cm of result, fill=green!10] 
        {Final verified response};

    \draw [arrow_style, dashed] (result.south) -- (final_output.north);
    
    \draw [arrow_style, blue, dashed, overlay] (result.east) to[out=0, in=-45] 
        node[midway, right, align=center, fill=white] {Neuro-symbolic\\feedback loop} 
        ($(stackbox.south east)!0.2!(stackbox.north east)$);

\end{tikzpicture}

\caption{Unlike standard LLMs which predict tokens directly from the transformer stack, an LRM generates a computational plan (Python script), executes it in an external deterministic sandbox, then integrates the result back into the neural stream for the final response. Verification is built-in. The model no longer predicts mathematical answers. It predicts the algorithms needed to find the answers and it executes the algorithms. }\label{LRM}
\end{figure}

\subsection{Structure of the point cloud of token embeddings in an LLM}
Let
\[
T^{(0)} := \{\varphi(w)\mid w \ \text{a token}\}\subset\mathbb{R}^{d}
\]
denote the  \emph{finite point cloud} of embedded tokens. 
The {\it Manifold Hypothesis of Machine Learning} \cite{TSL2000} asserts the existence of a smooth $n$-dimensional embedded
submanifold
\[
M^n\subset \mathbb{R}^{d}
\]
such that the probability distribution supported on $T^{(0)}$ sits on a 
neighborhood of $M$.  

In \cite{MDC2025}, the authors find that points in $T^{(0)}$ that are close in Euclidean distance exhibit `incompatible apparent dimensions'. That is, one point may have many nearby neighbors in $T^{(0)}$, while the other has far fewer.

This cannot occur if the point cloud is sampled from (or  near)
a single smooth embedded manifold $M^n\subset\mathbb{R}^{d}$. 
They conclude that there exists no low-dimensional smooth manifold $M$ in
$\mathbb{R}^{d}$ whose neighborhoods model the local statistical structure
of the token embedding point cloud. That is, the Manifold Hypothesis is violated.

Each transformer layer in the LLM implements a map
\[
f_\ell:\mathbb{R}^{d}\to\mathbb{R}^{d}
\]
built from linear projections (including $T_q$, $T_k$, $T_v$), attention
mixing, and a piecewise-linear feed-forward network.

Define the point cloud at layer $\ell$ to be
\[
T^{(\ell)} := f_\ell\circ f_{\ell-1}\circ\cdots\circ f_1(T^{(0)}).
\]

The observed behavior is compatible  with the point clouds
$T^{(\ell)}$ sampling neighborhoods of a \emph{stratified subset}
\[
S = \bigsqcup_i S_i \subset \mathbb{R}^{d},
\]
where the strata $S_i$ have different dimensions and meet along singular
loci.  The singular points are those points where the space is not smooth or locally Euclidean.
Near such points, no Euclidean ball intersects $S$ in a way resembling an open
ball in $\mathbb{R}^n$ for any $n$.

The authors show that the maps $f_\ell$ preserve the types of
singularities  already present in
$T^{(0)}$.  The singularities correspond to polysemous words (words with multiple meanings). The  different semantic regions of the embedded point cloud corresponding to polysemous words are `pinched' together \cite{JGZ2020}. Near the singularities, the LLM is unstable and `unusually sensitive', leading to unexpected associations and incoherent reasoning. This may be one of several causes of model hallucinations.

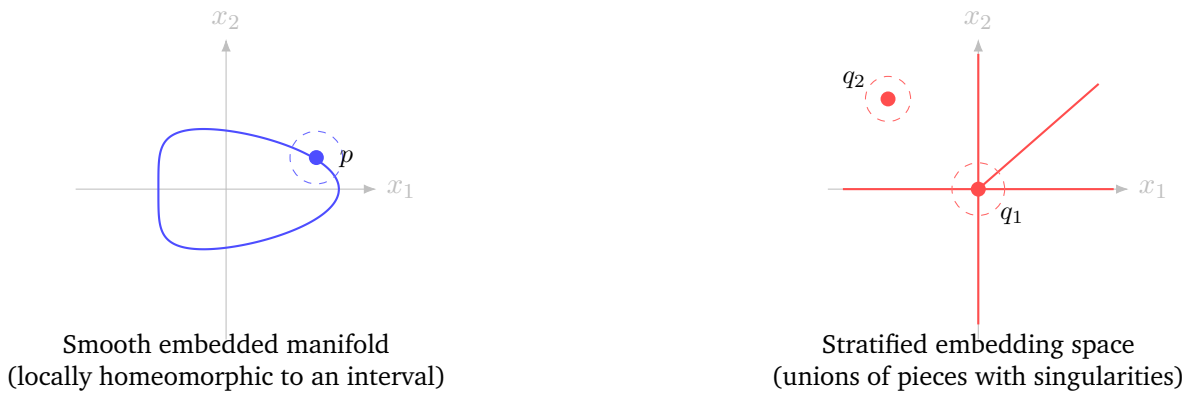
\begin{figure}[H]
\centering
\begin{tikzpicture}[   
  >=Latex,
    axis/.style={->, thin},
    manifold/.style={thick, blue!70},
    stratline/.style={thick, red!70},
    stratpoint/.style={circle, fill=red!70, inner sep=2pt},
    smoothpoint/.style={circle, fill=blue!70, inner sep=2pt},
    labelstyle/.style={font=\small, align=center}
]

\begin{scope}[shift={(-5,0)}]
    \draw[axis, gray!50] (-2,0) -- (2,0) node[right] {$x_1$};
    \draw[axis, gray!50] (0,-2) -- (0,2) node[above] {$x_2$};

    \draw[manifold, domain=0:360, smooth, samples=100]
        plot ({1.2*cos(\x) + 0.3*cos(2*\x)}, {0.8*sin(\x)});

    \node[smoothpoint] (p) at (1.2,0.42) {};
    \draw[thin, blue!60, dashed] (1.2,0.42) circle (0.35);
    \node[labelstyle] at (1.6,0.4) {$p$};

    \node[labelstyle] at (0,-2.3)
        {Smooth embedded manifold \\ (locally homeomorphic to an interval)};
\end{scope}

    {vs};

\begin{scope}[shift={(5,0)}]
    \draw[axis, gray!50] (-2,0) -- (2,0) node[right] {$x_1$};
    \draw[axis, gray!50] (0,-2) -- (0,2) node[above] {$x_2$};

    \draw[stratline] (-1.8,0) -- (1.8,0);
    \draw[stratline] (0,-1.8) -- (0,1.8);

    \draw[stratline] (0,0) -- (1.6,1.4);

    \node[stratpoint] (q) at (-1.2,1.2) {};

    \node[stratpoint] (s) at (0,0) {};

    \draw[thin, red!60, dashed] (0,0) circle (0.35);
    \draw[thin, red!60, dashed] (-1.2,1.2) circle (0.3);

    \node[labelstyle] at (0.45,-0.35) {$q_1$};
    \node[labelstyle] at (-1.65,1.45) {$q_2$};

    \node[labelstyle] at (0,-2.3)
        {Stratified embedding space \\ (unions of pieces with singularities)};
\end{scope}

\end{tikzpicture}
\caption{On the left,  (in blue) a smooth 1-dimensional manifold embedded in $\mathbb{R}^2$, where every point has a neighborhood homeomorphic to an open interval. On the right,  (in red) a stratified subset of $\mathbb{R}^2$, hypothetically resembling an LLM token embedding space. It is a union of 1-dimensional strata crossing at a singular point and an isolated 0-dimensional stratum. }
\end{figure}

\subsection{A single chat is a `context window'}

A generative AI model's short-term memory is its context window, which consists of a single chat session. This window contains the entire conversation history, including uploaded documents and any corrections made by the user. Its capacity is measured in tokens: most standard models offer a 120K window which can hold a novel. Models like Gemini  offer an unprecedented 1 million token input limit.

However, this memory is finite. Once the token limit is reached, the  AI model will begin to forget the earliest messages in the conversation. Work on a difficult mathematics problem must be contained within a single chat: the  AI model needs all information and previous steps to be present simultaneously to solve a complex problem.

Since a generative AI model  has no inherent memory between separate sessions, new chats must be initiated with any updated information and results.

\subsection{Large Reasoning Models and Large Context Models}
\label{LRM-LCM}
LLMs are predictive and fast. 
They `guess' the end result but they may not do all the intermediate work. 

{\it Large Reasoning Models (LRMs)}, are a class of LLMs  trained to carry out tasks that require multi-step  deduction, critical thinking, and structured problem-solving.
LRMs deliberate, explore various solution paths, evaluate intermediate steps and revise their reasoning. 

LRM models integrate frameworks that support reasoning structures, such as `Chain-of-Thought'  and `Tree-of-Thought' structures into their systems. We give some examples in Section~\ref{prompts}.

The training of LRMs diverges from that of LLMs.  They utilize  a technique called `process supervision'.   Unlike the training process of LLMs that only rewards a correct final answer, process supervision rewards the correctness of the intermediate reasoning steps.

LRM models also have self-correcting capabilities. Their basis is a neuro-symbolic system which combines text generation with verification (see Figure~\ref{LRM}). They have a feature known as `self-awareness'.  LRMs are able to distinguish between requests that are language-based or computational. They make `decisions' about delegating tasks to external tools such as a  Python interpreter. They don't  try to predict mathematical answers, they predict the algorithms needed to find the answers, then they execute the algorithms and report back with the result.

{\it Large Context Models (LCMs)} are defined by their  large context windows and the volume of information that they can handle. This allows them to accept and analyze large datasets, entire codebases, or  multimodal inputs, without  requiring external retrieval systems.
 Their training data  includes long-form content to improve their ability to track long-range dependencies between input tokens.

These advanced capabilities are becoming integrated into leading AI platforms. The cost, however,  is additional computing power and increased time for responses. 

Google's Gemini models demonstrate both LRM and LCM characteristics.  These models also use internal `thinking processes' for improved reasoning.

ChatGPT has a `reasoning mode' in its most recent model. Reinforcement learning was used to teach the model to `think' before generating answers, using what OpenAI refers to as  `private chain-of-thought' processes built into their system prompts. This allows the model to plan ahead and reason through tasks, performing a series of intermediate reasoning steps to assist in solving the problem. 

Claude models  use `extended thinking' or the ability to `think out loud'. 

These architectural advances are complemented by the practical integration of LLMs with symbolic computation tools.

\section{Computational Capabilities}\label{capabilities}
\subsection{AI models simulate Computer Algebra Systems}\label{CAS}
Advanced generative AI models simulate a scientific computing environment by delegating the computation to specialized Python libraries.

\begin{table}[h!]
\centering
\caption{Key Python libraries used by generative AI models for mathematics}
\label{tab:python_libs}
\begin{tabular}{@{} l l p{7.5cm} @{}}
\toprule
\textbf{Library} & \textbf{Main purpose} & \textbf{Key capabilities \& LLM use } \\
\midrule

\textbf{NumPy} & Numerical Computing & Handles arrays \& matrices. Used by generative AI models for  linear algebra (matrix multiplication, row reduction, eigenvalues). \\
\addlinespace

\textbf{SymPy} & Symbolic Mathematics & Performs precise algebra \& calculus. Used for symbolic derivatives, integrals, and simplifying expressions. \\
\addlinespace

\textbf{SciPy} & Scientific Computing & Provides advanced numerical routines. Used for optimization (finding minima/maxima) and solving differential equations. \\
\addlinespace

\textbf{Pandas} & Data Analysis & Manages structured, table-like data. Used for reading and analyzing data from files such as Excel. \\
\addlinespace

\textbf{Matplotlib} & Plotting \& Visualization & Creates a wide variety of  2D graphs and charts. Used to plot functions and visualize data. \\
\bottomrule
\end{tabular}
\end{table}
The Python code is written  by the LLM using  probabilistic methods. The code is then executed by a standard Python interpreter. The reliability of the final answer comes from having outsourced the computation to the verified, deterministic Python libraries. The potential unreliability comes from the LLM's probabilistic process of writing the code that uses those libraries.

For example, the LLM may hallucinate  functions that don't exist in the Python libraries. In addition, the queries may give incorrect answers if the problems are inaccurately posed. The LLM  could misunderstand the prompt, have bugs in its code or fail to handle boundary cases. Debugging AI-generated Python code may be needed.

Gemini and the paid version of ChatGPT   give automatic user access to the Python packages NumPy, SymPy, SciPy, Pandas, Matplotlib.

\subsection{Comparison of models}

Gemini, ChatGPT  and Claude  all exhibit  different specializations in how they activate nodes, navigate paths and  narrow down possibilities within their internal graphs of statistical associations.

For intermediate mathematical tasks, Gemini, ChatGPT, and Claude  are all roughly in the same  class with respect to performance. 
Though Gemini models performed  slightly better in  the standard benchmarks for competition-level mathematics, testing symbolic reasoning in algebra, calculus  and number theory in 2024.

ChatGPT integrates features like WolframAlpha, which is a powerful Computer Algebra System,  which can be enabled in the professional version.  However, WolframAlpha is not designed for multi-step algorithms and it times-out on long tasks. Also,  WolframAlpha is not optimized for integers. It   makes errors with integers larger than $10^{15}$. 

On the other hand, Python interpreters, available in Gemini and ChatGPT (see Subsection~\ref{CAS}), have `arbitrary-precision' for integers. They automatically adjust the memory used to store an integer as its value grows. We note that Python is not enabled in WolframAlpha.

Another ChatGPT feature called SciSpace has access to hundreds of millions of  peer-reviewed journal papers.

\section{Prompt engineering for mathematics}\label{prompts}
\subsection{How the prompt shapes an LLM's output}
The prompt provides an initial sequence of tokens that the model uses to condition its output. The LLM generates text by sequentially examining the existing sequence of tokens, assigning a probability to every possible next token, and then choosing one to add to the sequence.

The probability of each potential next token is calculated based on the full sequence of tokens that came before it -- both the original prompt and any text already generated. 

The final output is  the result of a chain of these probabilistic decisions. The framing of the initial prompt is critical in determining the usefulness of the output.

The design of prompts, or `prompt engineering', significantly influences the  behavior of an LLM \cite{PromptGuide}, \cite{PromptEngineering}. It disambiguates the input, leading to more relevant and useful output.  

Structured prompts impose a logical framework on an LLM's probabilistic generation process. The examples of  prompts in the following subsections have been shown to be useful in  guiding  generative AI models on how to approach complex mathematical questions \cite{kojima2022}. These methods are also used internally in advanced LRM models like Gemini and ChatGPT (in thinking mode).

\newpage
\subsection{$n$-shot prompts}

An {\it $n$-shot prompt} provides the LLM with $n$ examples of the task that you would like it to perform.

\bigskip
 {\bf{Example: Zero-shot prompt}}

{\color{blue} Determine the arithmetic progression corresponding to primes of the form $4n+1$ with last digit $1$,  (such as $41, 61, 101, \dots$).}

\bigskip
 {\bf{Example: One-shot prompt}}

{\color{blue} Given that primes of the form $4n+1$ ending in $1$ (such as $41, 61, 101, \dots$) correspond to the arithmetic progression ${20k + 1}$, determine the arithmetic progression corresponding to primes of the form $4n+1$ with last digit $3$ (such as $13, 53, 73, 113, \dots$).}

Giving the LLM an example makes it more likely to give a correct answer.

\subsection{Chain-of-Thought reasoning}

This refers to the  LLM thinking step-by-step, generating intermediate reasoning before arriving at the final answer  \cite{IBMToT}.  If the LLM checks  its work, this reduces errors.

{\bf Example prompt:}

{\color{blue} Your task is to generate a detailed description of the representation of the Weyl group of $\mathfrak{sl}_3(\mathbb{C})$ on the dual  space $\mathfrak{h}^*$ of its Cartan subalgebra $\mathfrak{h}$. 

Start by describing the Coxeter presentation of the group $W(A_2) \cong S_3$ which has generators $s_1, s_2$ and relations $s_i^2=e$ and $(s_1s_2)^3=e$. Next, define the vector space $\mathfrak{h}^*$ and specify its basis of simple roots, $\alpha_1, \alpha_2$. 

Then define the representation $\rho$ of $W(A_2)$ using the fact that the generators act as  reflections on $\mathfrak{h}^*$. 

Use the Weyl reflection formula to explicitly derive the matrices for $\rho(s_1)$ and $\rho(s_2)$ in the basis for $\mathfrak{h}^*$. 

With this information, the next step is to verify that this is a valid representation. Please do this by performing the necessary matrix calculations to show that these matrices satisfy the  defining relations of the group $W(A_2) \cong S_3$. 
}

An LRM may carry out internal reasoning as in Figure~\ref{CoT}. A LLM can be explicitly prompted to use Chain-of-Thought reasoning.

\begin{figure}
\begin{center}
\begin{tikzpicture}[
  node distance=10mm and 12mm,
  >=Latex,
  block/.style={rounded corners, draw, very thick, align=center, inner sep=6pt, fill=gray!5},
  title/.style={draw, very thick, align=center, inner sep=6pt, fill=blue!5},
  smallmath/.style={font=\small}
]
\node[title] (goal) {{\bf Goal}\\ Representation of \(W(A_2)\cong S_3\) on \(\mathfrak{h}^*\subset\mathfrak{sl}_2(\C)\) };
\node[block, below=of goal] (coxeter) {\textbf{Weyl group}\\[2pt]\(\displaystyle W(A_2)=\langle s_1,s_2\mid s_i^2=e,\ (s_1s_2)^3=e\rangle\)};
\node[block, below=of coxeter] (hstar) {\textbf{Basis for \(\mathfrak{h}^*\)}\\[1.5pt]Specify basis of simple roots \(\{\alpha_1,\alpha_2\}\),\ Cartan matrix \(A=\begin{pmatrix}2&-1\\-1&2\end{pmatrix}\)};
\node[block, below=of hstar] (action) {\textbf{Define action by reflections}\\[2pt]For \(\lambda\in\mathfrak{h}^*\), \(\displaystyle s_i(\lambda)=\lambda-\langle \lambda,\alpha_i^\vee\rangle\,\alpha_i.\)\\on simple roots: \(\displaystyle s_i(\alpha_j)=\alpha_j-a_{ij}\alpha_i\)};
\node[block, below=of action] (matrices) {\textbf{Compute matrices in basis \((\alpha_1,\alpha_2)\)}\\[4pt]\(\displaystyle \rho(s_1)=\begin{pmatrix}-1 & 1\\[2pt]0 & 1\end{pmatrix},\qquad \rho(s_2)=\begin{pmatrix}1 & 0\\[2pt]1 & -1\end{pmatrix}.\)};
\node[block, below=of matrices] (verify) {\textbf{Verify relations}\\[2pt]\(\rho(s_i)^2=I_2\)\quad and \quad\((\rho(s_1)\rho(s_2))^3=I_2\)};
\node[title, below=of verify] (conclude) {{\bf Conclusion}\\ \(\rho:W(A_2)\to \mathrm{GL}(\mathfrak{h}^*)\) is a representation};
\draw[->, very thick] (goal) -- (coxeter);
\draw[->, very thick] (coxeter) -- (hstar);
\draw[->, very thick] (hstar) -- (action);
\draw[->, very thick] (action) -- (matrices);
\draw[->, very thick] (matrices) -- (verify);
\draw[->, very thick] (verify) -- (conclude);
\node[block, smallmath, right=20mm of action] (rule) {Use \(A=(a_{ij})\) so that\\\(s_i(\alpha_j)=\alpha_j-a_{ij}\alpha_i\).};
\draw[->, thick] (action.east) -- (rule.west);
\node[block, smallmath, right=20mm of matrices, yshift=-2mm] (prod) {Compute \(\rho(s_1)\rho(s_2)\)\\and check \((\rho(s_1)\rho(s_2))^3=I_2\).};
\draw[->, thick] (matrices.east) -- (prod.west);
\end{tikzpicture}
\end{center} \caption{Chain-of-thought diagram.}\label{CoT}
\end{figure}
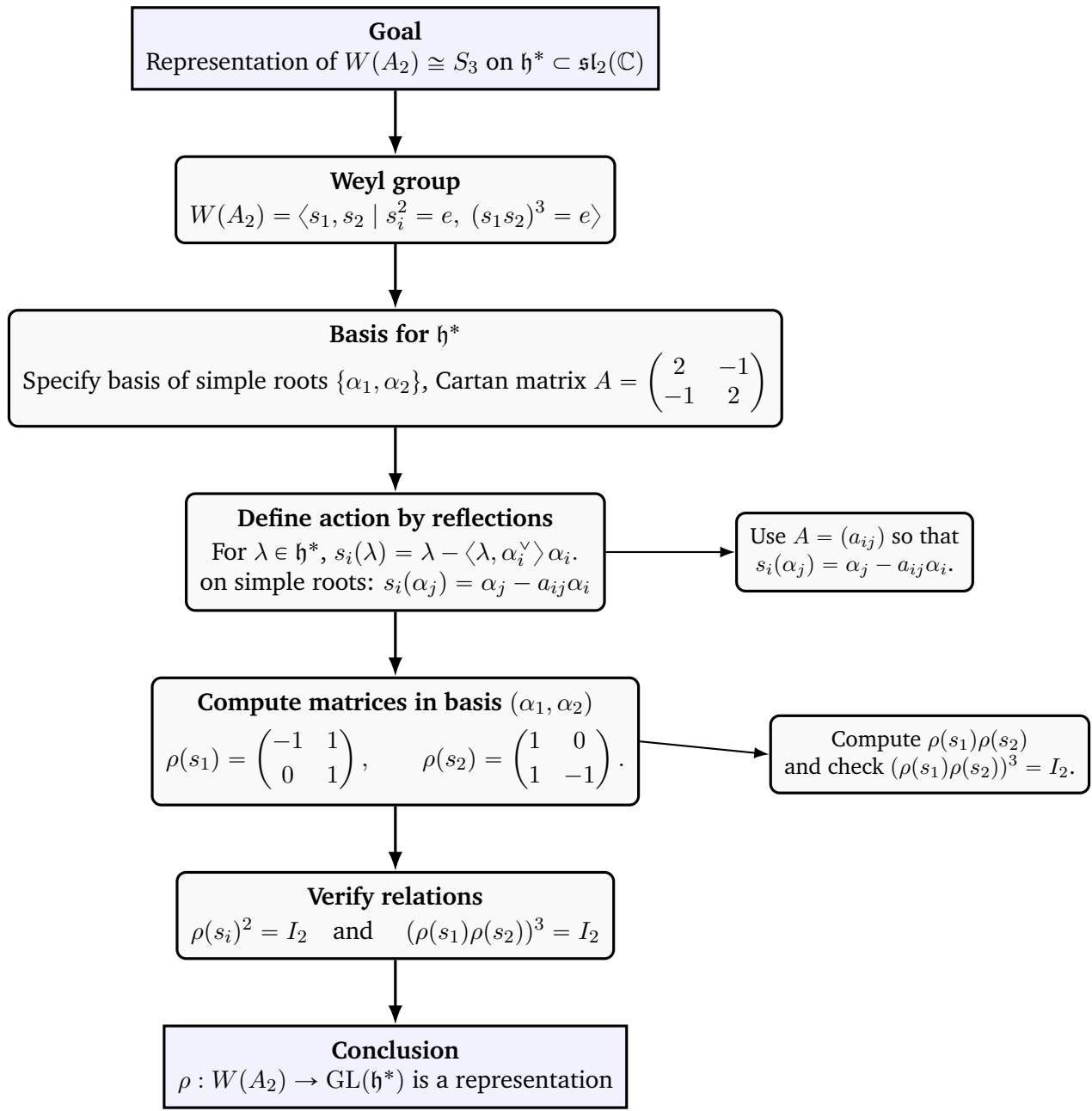

\clearpage

\subsection{Tree-of-Thought reasoning}

This refers to an AI model exploring multiple different deduction paths simultaneously, such that it evaluates its own progress at each step  and pursues the most promising path   \cite{IBMToT}.

{\bf Example prompt:}

{\color{blue} Let $p$ be a prime number, written $p=a_na_{n-1}\dots a_1a_0$ in terms of its digits in base 10. Suppose that $n>1$, that is, suppose that $p$ has at least two digits. Your task is to determine the possible values of the last digit $a_0$. }

An LRM may carry out internal reasoning as in Figure~\ref{ToT}. A LLM can be explicitly prompted to use Tree-of-Thought reasoning.
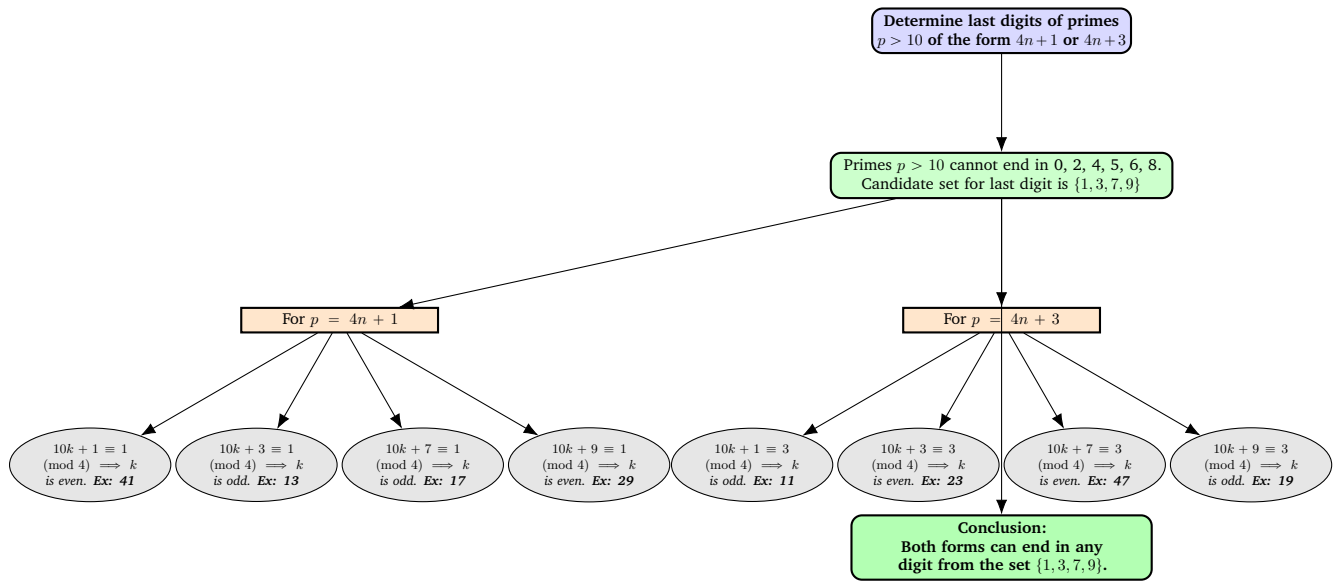
\begin{figure}[h!]
\centering
\begin{tikzpicture}[
    scale=0.55, transform shape,
    level distance=3.5cm,
    level 1/.style={sibling distance=0cm}, 
    level 2/.style={sibling distance=16cm},
    level 3/.style={sibling distance=4.02cm},
    root/.style={
        rectangle, rounded corners, draw, thick,
        fill=blue!15,
        text width=6cm, align=center, font=\bfseries
    },
    filter/.style={
        rectangle, rounded corners, draw, thick,
        fill=green!20,
        text width=8cm, align=center
    },
    branch/.style={
        rectangle, draw, thick,
        fill=orange!20,
        text width=4.5cm, align=center
    },
    leaf/.style={
        ellipse, draw,
        fill=gray!20,
        text width=2.5cm,
        align=center, font=\footnotesize\itshape
    },
    final/.style={
        rectangle, rounded corners, draw, thick,
        fill=green!30,
        text width=7cm, align=center, font=\bfseries
    },
    edge from parent/.style={draw, -{Latex[length=2mm]}}
]
\node[root] {Determine last digits of primes $p >10$ of the form $4n+1$ or $4n+3$}
    child {
        node[filter] {Primes $p > 10$ cannot end in 0, 2, 4, 5, 6, 8. \\ Candidate set for last digit is $\{1, 3, 7, 9\}$}
        child {
            node[branch] {For $p = 4n+1$}
            child { node[leaf] {$10k+1 \equiv 1  \pmod 4 \implies k$ is even. \textbf{Ex: 41}} }
            child { node[leaf] {$10k+3 \equiv 1 \pmod 4  \implies k$ is odd. \textbf{Ex: 13}} }
            child { node[leaf] {$10k+7 \equiv 1 \pmod 4  \implies k$ is odd. \textbf{Ex: 17}} }
            child { node[leaf] {$10k+9 \equiv 1 \pmod 4  \implies k$ is even. \textbf{Ex: 29}} }
        }
        child {
            node[branch] {For $p = 4n+3$}
            child { node[leaf] {$10k+1 \equiv 3 \pmod 4 \implies k$ is odd. \textbf{Ex: 11}} }
            child { node[leaf] {$10k+3 \equiv 3 \pmod 4 \implies k$ is even. \textbf{Ex: 23}} }
            child { node[leaf] {$10k+7 \equiv 3 \pmod 4 \implies k$ is even. \textbf{Ex: 47}} }
            child { node[leaf] {$10k+9 \equiv 3 \pmod 4 \implies k$ is odd. \textbf{Ex: 19}} }
        }
        child[level distance=9cm, sibling distance=0cm] {
            node[final] {Conclusion:\\ Both forms can end in any digit from the set $\{1, 3, 7, 9\}$.}
        }
    };
\end{tikzpicture}
\caption{A Tree-of-Thought diagram for determining the possible last digits of primes $p>10$.}\label{ToT}
\end{figure}

\subsection{Using LLMs to generate Python code to run in SageMath}

Generative AI models can generate code for Computer Algebra Systems such as GAP, Magma, and SageMath with varying degrees of reliability. Proficiency is particularly high with SageMath, given its Python-based syntax and integration with GAP for computations in group theory, ring theory and field theory.

{\bf Example prompt:}

{\color{blue} Generate Python code to find the minimum number of colors needed to color the Petersen Graph such that no two adjacent vertices share the same color (the `chromatic number'). Then create the Petersen Graph using the built-in SageMath generator. Find its chromatic number. Display the graph.
}

\bigskip
\begin{figure}[H]\begin{center}
	\includegraphics[width=3.2in]{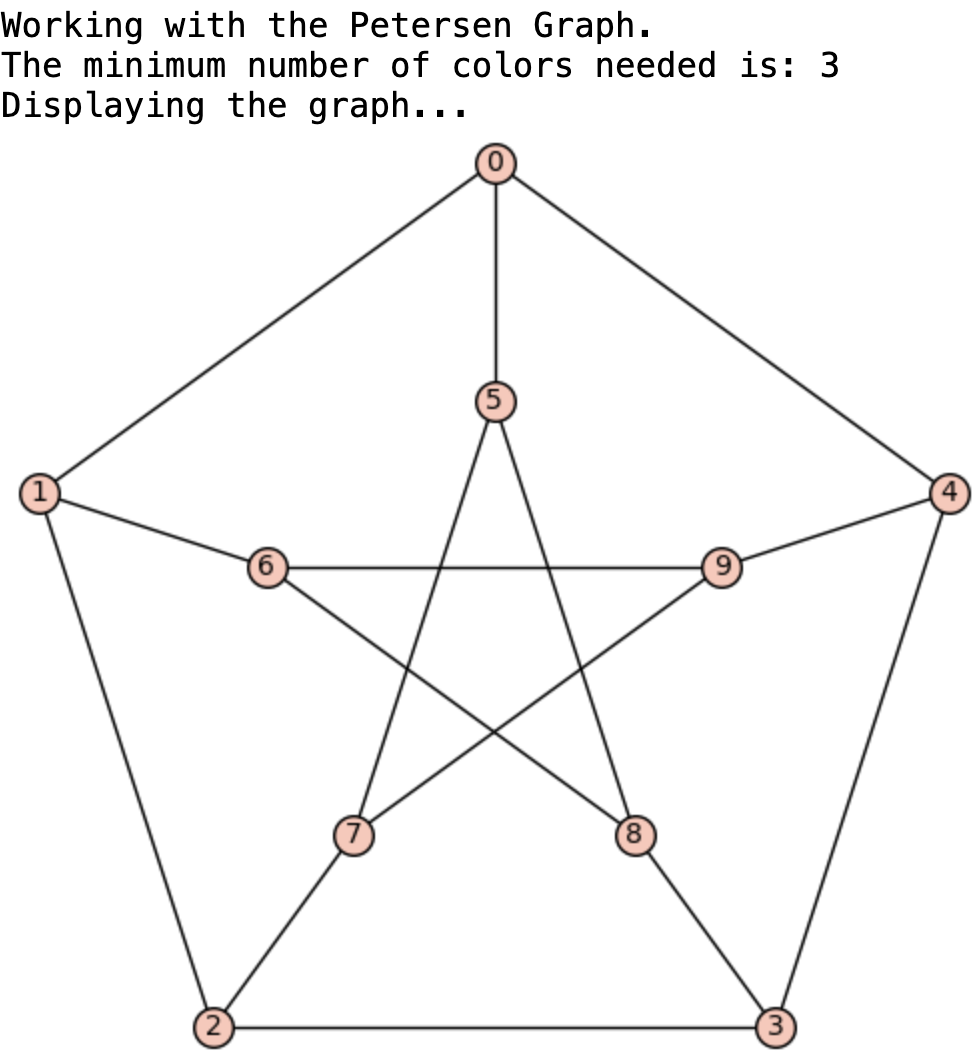}
		\end{center}
	\end{figure}

The emergence of LRMs and LCMs (Section~\ref{LRM-LCM}) has  altered the practice of prompt engineering.

For LRMs, the need for  prompt engineering techniques designed to force reasoning, such as explicit Chain-of-Thought prompting, is reduced. The reasoning process is handled internally. Simpler, goal-oriented prompts that clearly define the problem may be more effective.

With LCMs, the focus is shifting to `context engineering'. The challenge becomes optimizing the content and structure of the  information provided in the prompt. 
This  maximizes performance and prevents recall accuracy from decreasing as the context window fills.

\subsection{Targeted prompts}
You can ask a generative AI model to:

\begin{itemize}
\item Give an explicit reference (with page numbers) on a topic known to be in some book or research paper.

\item Apply the statement of a theorem to a particular example.

\item Generate LaTeX code from an uploaded pdf file of a math paper. 

\item Translate a mathematics paper from another language  and give the output in a LaTeX file.

\item Explain a section of a physics paper in mathematical language. 

\item Generate a Tikz diagram or table from a description in natural language.

\item Generate a bibliography on a specific topic.

\item Rewrite the LaTeX code for an entire math paper in different notation.

\item Find the typos in an uploaded pdf file.

\item Analyze if the flow of ideas in a paper is appropriately sequential.

\item Generate LaTex, Lean, Python, Mathematica, Maple and other forms of code.

\end{itemize}

\newpage
\subsection{General rules for writing prompts for advanced math questions}
\begin{itemize}

\item Specify the task in short sentences. 

\item Use `your task is to...' or  `your goal is to...'. 

\item Specify the context of the task.

\item Upload all necessary background information.

\item Be explicit and detailed.

\item Include all relevant keywords.

\item Guide the reasoning process.

\item Specify the output format.

\item Verify the output.

\item Never ask for a complex proof in one shot.

\item Break a difficult task up into smaller tasks.

\item Take the model's output and ask for it to be modified  with specific constraints.

\item Constrain the method the model is allowed to use.

\item  Find errors in the output and ask the model to self-correct. 

\end{itemize}
 You can also ask a generative AI model to generate a prompt for a given task.

\section{Other ways to influence the  output from a generative AI model}\label{output}

\subsection{Design your outputs}
In generative AI models, there is a  `system instruction' feature which can be used to influence the  output. Users can enter a `pre-prompt' with their own personal profile that governs the model's behavior for the entire chat session. Some models allow users to build AI agents that have persistent memory of settings. This capability allows for a type of `back-end' engineering, where writing specific system prompts gives the user more direct control over the model's performance. Utilizing this feature can significantly impact the amount of detail provided in the response and the level of rigor maintained in a mathematical proof.


\subsection{Changing the `knobs'}  System settings in certain generative AI models can also be manually changed in order to make mathematical proofs more rigorous and less random.

For mathematical rigor, the most important setting is  `temperature'. This
controls the randomness of the output by reshaping the probability distribution. For research questions and formulating conjectures, a high temperature (such as 0.9) is preferable, as it encourages more diverse and novel responses. A low temperature (such as 0.2 or 0.1) makes the output more deterministic. 

\begin{table}[h!]
\centering
\caption{Optimized settings in Google AI studio for research in mathematics}
\label{tab:settings}
\resizebox{\textwidth}{!}{%
\begin{tabular}{@{}l p{0.8\textwidth}@{}}
\toprule

\textbf{Temperature} & \textbf{Medium-High} ($0.7$ - $1.0$) \par  A higher temperature encourages exploration. 
\\

\addlinespace

\textbf{Top{ }P} & \textbf{High} ($0.95$ - $1.0$) \par  Allows the model to consider a wider, more diverse set of `next words'. 
\\
\addlinespace

\textbf{Thinking mode} & \textbf{Advanced} \par  Research problems  require the deepest level of reasoning. 
\\
\addlinespace

\textbf{Set thinking budget} & \textbf{Maximum} \par  Gives the model time to explore the `search space'. 
\\
\addlinespace

\textbf{Code execution} & \textbf{On} \par  Crucial for running checks with Python tools such as SymPy. 
\\
\addlinespace

\textbf{Grounding (Google Search)} & \textbf{On} \par  Allows the model to search outside its training data.
\\
\addlinespace

\textbf{URL Context} & \textbf{On} \par   Allows you to give the model specific URLs as context prompts.\\
\addlinespace

\textbf{Structured output} & \textbf{Off} \par  Allows unstructured brainstorming. 
\\
\addlinespace

\textbf{Function calling} & \textbf{On} \par  Connect to external research tools such as a university library API or a computational algebra system like SageMath or WolframAlpha. 
\\
\addlinespace

\textbf{Media resolution} & \textbf{High} \par  Essential for correct reading of symbols in handwritten notes or complicated diagrams. 
\\
\addlinespace

\textbf{Output Length} & \textbf{Maximum} \par  Allows for an unencumbered output. 
\\
\bottomrule
\end{tabular}
}
\end{table}

`Back-end' engineering, such as writing  system prompts in a generative AI model gives more direct control and can impact the amount of detail and rigor in a proof.

\subsection{Fact-checking: retrieval-augmented generation} 

The default output of a generative AI model is an output that is based on its training data.
The ability to access the  internet to fact-check in real time is a feature built on top of an LLM,  through a mechanism called `retrieval-augmented generation', or RAG. It may require a specific prompt in order to invoke this feature.
Users can also explicitly request a list of websites and references used.

When asked a question that requires fact-checking, Gemini
automatically queries the Google Search index in real time. It then synthesizes the information from the  top-ranked pages to give an answer and often provides direct links.

The current version of ChatGPT has a `reasoning mode' in its new model which autonomously decides when to use its web-browsing tools. Before finalizing an answer, the model evaluates and corrects its own results. It  internally grades its answer against criteria to ensure it meets a high-quality standard.

Microsoft Copilot is designed  to be an AI-powered search assistant. Most queries  initiate a live search on Bing and links are usually provided.

\section{Case Study: Combinatorial Group Theory}\label{CGT}
\subsection{Why AI is a useful tool in Combinatorial Group Theory}

We discuss some examples from group theory. However,  many of the underlying principles are applicable to other topics.

Let $G$ be a group defined by a presentation $G=\langle X \mid R \rangle$, where $X$ is a set of generators and $R$ is a set of relations. This is a compact way to define a group, but it hides immense complexity. 

The word problem (determining if a word in the generators is equal to the trivial element of the group) is unsolvable in general \cite{MT1973}. 

Computer Algebra Systems such as GAP, Magma, and SageMath can  solve the word problem for specific families of groups using  rewriting methods and structural properties. 

A  simple example: let $G  = BS(1,2)= \langle x, y \mid x^{-1}yx = y^2 \rangle$. 
This is the so-called {\it Baumslag-Solitar group} $BS(1,2)$. Its word problem is known to be solvable.

A solution of the word problem involves deciding on an efficient  sequence of applications of the group relations to achieve the goal of trying to get to the identity element of the group.

As an example, we solve the word problem for the word $w = yx y^{-2} x^{-1}$:

We first rewrite $x^{-1}yx = y^2$ as $yx = xy^2$. Then 
$${yx} y^{-2} x^{-1}={xy^2}y^{-2} x^{-1}=xx^{-1}=1.$$
Thus the word $w = yx y^{-2} x^{-1}$ represents the identity element in $G$.

\subsection{Solving algorithmic research questions with AI}

Algorithmic search problems in group theory, such as the word problem, conjugacy problem (deciding if two words are conjugate in a given group), and triviality problem (deciding if a group presentation presents the trivial group), are undecidable in general. Even in decidable cases, the search space for groups with many generators and intertwined relations is often astronomically large. The number of ways to apply relations to a word grows exponentially with the length of the word.

AI can address this by effectively simplifying words. It has strategies to navigate the  search space by deciding which group relations are the most useful to apply. 

Gemini uses {Chain-of-Thought} processes for working with group presentations. The model does not just guess the final simplified word. It generates intermediate steps internally using a type of depth-first search where the model explores a simplification path. It can self-correct if the word complexity increases rather than decreases.

 Gemini translates the group presentation $G=\langle X\mid R\rangle$ into a Python script using the \texttt{sympy.combinatorics} module.
    \begin{itemize}
        \item It defines the free group on generators $X$.
        \item It inputs the relations $R$.
        \item It uses the library's implementation of the {Knuth-Bendix completion algorithm} (or similar rewriting systems) to attempt to find a normal form for the word.
    \end{itemize}

The following is an example of a prompt explicitly asking the model to follow  Tree-of-Thought reasoning.

{\bf Example prompt:}

{\color{blue}Your task is to determine if the Higman group, defined by the presentation 
\[
    H = \langle a, b, c, d \mid a^{-1}ba=b^2, \, b^{-1}cb=c^2, \, c^{-1}dc=d^2, \, d^{-1}ad=a^2 \rangle
\]
is trivial or non-trivial. Proceed with the following tree of reasoning, evaluating each step before proceeding to the next:}

\begin{figure}[H]
\centering
\begin{tikzpicture}[
    scale=0.69, transform shape,
    level distance=2.5cm,
    sibling distance=13cm, 
    level 2/.style={sibling distance=5.5cm}, 
    root/.style={
        rectangle, rounded corners, draw, thick,
        fill=blue!15,
        text width=4.5cm, align=center, font=\bfseries
    },
    branch/.style={
        rectangle, draw, thick,
        fill=orange!20,
        text width=3cm, align=center
    },
    leaf/.style={
        ellipse, draw,
        fill=gray!20,
        text width=2.5cm,
        align=center, font=\footnotesize\itshape
    },
    final/.style={
        ellipse, draw, thick,
        fill=green!30,
        text width=4.5cm, align=center, font=\bfseries\large
    },
    edge from parent/.style={draw, -{Latex[length=2mm]}}
]
\node[root] {Is the Higman Group Trivial?}
    child {
        node[branch] {Branch 1: Elementary Tests}
        child {
            node[branch] {1a. Abelianization}
            child { node[leaf] {$H_{ab}$ is trivial. No non-trivial abelian quotients.\\ \textbf{Inconclusive}} }
        }
        child {
            node[branch] {1b. Direct Simplification}
            child[level distance=3.5cm] { 
                node[leaf] {Leads to complexity. \\ \textbf{Inconclusive}} 
            }
        }
    }
    child {
        node[branch] {Branch 2: Structural Analysis}
        child {
            node[branch] {2a. Identify group generated by individual relations}
            child[level distance=3.5cm] { 
                node[leaf] {Identifies non-trivial $BS(1,2)$ building blocks.} 
            }
        }
        child {
            node[branch] {2b. Explain Construction}
            child { node[leaf] {Built via  HNN extensions with embedded $BS(1,2)$ subgroups.} }
        }
        child[level distance=3.5cm] {
            node[final] {Conclusion:\\ The Higman group is non-trivial.}
        }
    };
\end{tikzpicture}
\caption{Tree-of-Thought prompt diagram for non-triviality of the Higman group.}
\end{figure}
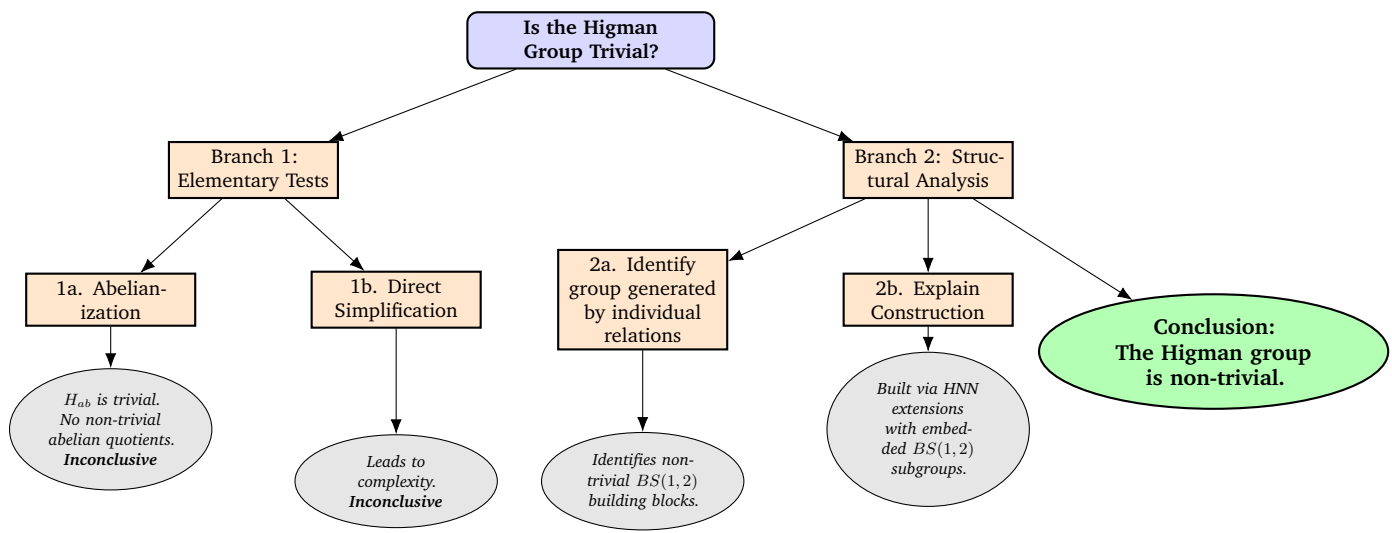

\section{Hybrid systems: generative AI models with Computer Algebra Systems or formal proof assistants}\label{hybrid}

Because generative AI models are not substitutes for formal proof or verified computation, their outputs require rigorous human scrutiny. The primary way to offset this limitation is to use them in conjunction with formal tools like Computer Algebra Systems or formal proof assistants such as Lean. For example, the professional version of ChatGPT can integrate with WolframAlpha. The user can explicitly prompt ChatGPT in the Wolfram chat window, or allow it to decide when to outsource a computation.

An active area of research is to develop tools that link generative AI models with Computer Algebra Systems \cite{KG2025}.

Formal proof assistants like Lean allow for the expression and mechanical verification of mathematical proofs. Lean implements a version of Dependent Type Theory known as the Calculus of Inductive Constructions (CIC).  This is a constructive type theory, where in the Lean kernel (or core verifier)  proofs are given explicitly, aligning with computation and the Curry--Howard correspondence. That is, every proof encodes an algorithm or a construction. However, in practice, Lean contains a model of ZFC set theory. It  does not use ZFC axioms to check proofs, but with some additional axioms, it behaves like ZFC.

Type Theory is particularly advantageous for formalization because it embeds mathematical meaning directly into its syntax.   Its system of Types prevents false statements and builds properties directly into the definition of its objects.  In ZFC set theory, many of these logical statements require separate proofs.

However, writing Lean code is notoriously difficult, requiring not just programming skill but the precise formalization of abstract mathematical ideas in Type Theory. In this framework, even small logical gaps or unstated assumptions must be made explicit and rigorously verified. Precise  use of syntax is 
required. 

Most users use Lean {\it tactics}. These  are metaprograms that automate steps and give instructions on how to construct proofs. Meanwhile,  Lean generates the rigorous low-level code for proof-checking.

One of the benefits of the Lean environment is  Mathlib: a vast, open-source digital encyclopedia of formally verified mathematics.

Lean Copilot is a framework that integrates LLMs  into the Lean theorem prover to assist  with formal proofs. It functions as an AI assistant by providing features like context-aware tactic suggestions. It selects relevant premises from Mathlib and searches for proofs.

Research labs are actively working on ways to make the task of  generating Lean code easier: from using generative AI models to write Lean code, to  developing interfaces that accept inputs in natural language. 

For example, {\it Numina} is an assistant that acts as a bridge between natural language mathematics proofs and  Lean. It uses AI to generate formal proofs and uses Lean to verify their correctness.

\section{Mathematician $+$ generative AI collaboration}

\subsection{Human-LLM collaborative research}
Using a generative AI model as a collaborative research partner, rather than as a search engine or text generator, involves an  iterative dialogue that must be repeatedly corrected and refined. The user can improve the model's  output through a series of increasingly specific constraints, corrections, and questions, often in a loop that involves downloading, correcting, and re-uploading the model's work for further prompting.

For open research questions, this can be a lengthy process, sometimes involving hundreds of prompts  in a single chat, to obtain a final product that is useful. 
A key technique is `in-context learning', where researchers provide background information like papers or books via prompts. The LLM keeps this material in its short-term memory. It  becomes more adept at handling the research topic as the conversation history forms an expanding context window.

While not necessarily a time-saving endeavor, this  collaborative process can reveal ideas and connections that lie beyond a human researcher's own spectrum.

\subsection{Why mathematicians should collaborate with generative AI models}
 As we have discussed, generative AI, on its own, should not be viewed as an oracle or an authoritative source for mathematical output. Verification of output  is a necessary part of using an AI model for mathematics.  But it can be used as a creative collaborator that provides ideas, perspective and potential new directions. 
 
 The most effective research strategy uses all available tools. When combined with a Computer Algebra System, AI can reliably perform laborious tasks and calculations. When paired with formal verification tools, it transforms the research process itself. 
 
 The author of this work now employs this collaborative framework in  several LLM+CAS+Lean-assisted research projects in infinite dimensional algebra and group theory. The use of generative AI models has been particularly revealing in this context. There have been multiple instances where AI models have suggested  novel research directions and ideas that  had initially been overlooked or dismissed as implausible or irrelevant. Yet, these suggestions (once verified) turned out to yield   unexpected paths forward. 

Mathematicians could consider AI not merely as a search engine, but as a creative partner that can act as a catalyst for discovery.

\newpage
\bibliography{AIpaper}{}

\end{document}